\documentclass[a4paper,12pt]{article}
\usepackage{latexsym}
\usepackage{amsfonts}
\newtheorem{theorem}{Theorem}[section]
\newtheorem{corollary}[theorem]{Corollary}

\newtheorem{example}[theorem]{Example}
\newtheorem{remark}[theorem]{Remark}
\newtheorem{lemma}[theorem]{Lemma}
\newtheorem{proposition}[theorem]{Proposition}
\title{On a stronger form of hereditary compactness in product
spaces\thanks{1991 Math.\ Subject Classification --- Primary:
54B10, 54D30; Secondary: 54A05, 54G99. \protect\newline Key words and
phrases --- sg-compact, hereditarily compact, $C_2$-space,
semi-open, sg-open, sg-closed, hsg-closed.}}
\author{Julian Dontchev\\Department of Mathematics\\University
of Helsinki\\PL 4, Yliopistonkatu 15\\00014 Helsinki 10\\Finland
\and Maximilian Ganster\\Department of Mathematics\\Graz
University of Technology\\Steyrergasse 30\\A-8010 Graz\\Austria}
\date{}
\begin{document}
\baselineskip=20pt plus 1pt minus 1pt
\newcommand{\fxy}{$f \colon (X,\tau) \rightarrow (Y,\sigma)$}
\maketitle
\begin{abstract}
The aim of this paper is to continue the study of sg-compact
spaces. The class of sg-compact spaces is a proper subclass of
the class of hereditarily compact spaces. In our paper we shall
consider sg-compactness in product spaces. Our main result says
that if a product space is sg-compact, then either all factor
spaces are finite, or exactly one factor space is infinite and
sg-compact and the remaining ones are finite and locally
indiscrete.
\end{abstract}

\section{Introduction}\label{s1}

If a topological space $(X,\tau)$ is hereditarily compact, then
under some additional assumptions either $X$ or $\tau$ might
become finite (or countable). For example, if $(X,\tau)$ is a
second countable hereditarily compact space, then $\tau$ is
finite. Hence, if $(X,\tau)$ is a second countable hereditarily
compact $T_0$-space, then $X$ must be countable. Moreover, it is
well-known that every maximally hereditarily compact space and
every hereditarily compact Hausdorff (even kc-) space is finite.
For more information about hereditarily compact spaces we refer
the reader to A.H. Stone's paper \cite{STO1}.

In 1995 and in 1996, a stronger form of hereditary compactness
was introduced independently in three different papers. Caldas
\cite{CC1}, Devi, Balachandran and Maki \cite{DBM1} and Tapi,
Thakur and Sonwalkar \cite{T2} considered topological spaces in
which every cover by sg-open sets has a finite subcover. These
spaces have been called {\em sg-compact} and were further studied
by the present authors in \cite{DG1}.

As the property sg-compactness is much stronger than hereditary
compactness (for even spaces with finite topologies need not be
sg-compact), the general behavior of sg-compactness becomes more
`unusual' than the one of hereditarily compact spaces. This will
be especially the case in product spaces.

It is well-known that the finite product of hereditarily compact
spaces is hereditarily compact, and that if a product space
is hereditarily compact, then every factor space is hereditarily
compact. What we want to show here is the following: If the
product space of an arbitrary family of spaces is sg-compact,
then all but one factor spaces must be finite and the remaining
one must be (at most) sg-compact. Maki, Balachandran and Devi
\cite[Theorem 3,7]{MBD1} showed (under the additional assumption
that the product space satisfies the weak separation axiom
$T_{gs}$) that if the product of two spaces is sg-compact, then
every factor space is sg-compact. Tapi, Thakur and Sonwalkar
\cite[Theorem 2.7]{T2} stated the result for two spaces but their
proof is wrong as they claimed that the projection mapping is
sg-irresolute. They used a wrong lemma from \cite{T1} saying that
the product of sg-closed sets is sg-closed (we will show that
this is not true even for two sets).

We recall some definitions. A set $A$ is called {\em semi-open}
if $A \subseteq {\rm cl} ({\rm int} (A))$ and {\em semi-closed}
if ${\rm int} ({\rm cl} (A)) \subseteq A$. The {\em
semi-interior} (resp.\ {\em semi-kernel}\/) of $A$, denoted by
${\rm sint}(A)$ (resp.\ ${\rm sker} (A)$, is the union (resp.\
intersection) of all semi-open subsets (resp.\ supersets) of $A$.
The {\em semi-closure} of $A$, denoted by ${\rm scl}(A)$, is the
intersection of all semi-closed supersets of $A$. A subset $A$
of a topological space $(X,\tau)$ is called {\em sg-open}
\cite{BL1} (resp.\ {\em g-open} \cite{L1}) if every semi-closed
(resp.\ closed) subset of $A$ is included in the semi-interior
(resp.\ interior) of $A$. A topological space $(X,\tau)$ is
called {\em sg-compact} \cite{CC1,DBM1,T2} (resp.\ {\em
go-compact} \cite{BSM1}) if every cover of $X$ by sg-open (resp.\
g-open) sets has a finite subcover.

Complements of sg-open sets are called {\em sg-closed}.
Alternatively, a subset $A$ of a topological space $(X,\tau)$ is
called sg-closed if ${\rm scl} (A) \subseteq {\rm sker} (A)$. If
every subset of $A$ is also sg-closed in $(X,\tau)$, then $A$ is
called {\em hereditarily sg-closed} (= hsg-closed) \cite{DG1}.
Every nowhere dense subset is hsg-closed but not conversely.

Jankovi\'{c} and Reilly \cite[Lemma 2]{JR1} pointed out that in
an arbitrary topological space every singleton is either nowhere
dense or locally dense. Recall that a set $A$ is said to be {\em
locally dense} \cite{CM1} (= {\em preopen}) if $A \subseteq {\rm
int} ({\rm cl} (A))$. We will make significant use of their
result throughout this paper.

\begin{lemma}\label{l0}
For a topological space $(X,\tau)$ the following conditions are
equivalent:

{\rm (i)} $X$ is locally indiscrete.

{\rm (ii)} Every singleton is locally dense.

{\rm (ii)} Every subset is sg-open.
\end{lemma}

\begin{lemma}\label{l1}
{\rm (i)} Every open continuous surjective function is
pre-semi-open, i.e., it preserves semi-open sets.

{\rm (ii)} Let $(X_i)_{i \in I}$ be a family of spaces and
$\emptyset \not= A_i \subseteq X_i$ for each $i \in I$. Then,
$\prod_{i \in I} A_i$ is preopen (resp.\ semi-open) in $\prod_{i
\in I} X_i$ if and only if $A_i$ is preopen (resp.\ semi-open)
in $X_i$ for each $i \in I$ and $A_i$ is non-dense (resp.\ $A_i
\not= X_i$) for only finitely many $i \in I$.

(iii) If \fxy\ is open and continuous, then
the preimage of every nowhere dense subset of $Y$ is nowhere
dense in $X$, i.e., $f$ is $\delta$-open.
\end{lemma}

\begin{lemma}\label{l2}
{\em \cite[Theorem 2.6]{DG1}} For a topological space $(X,\tau)$
the following conditions are equivalent:

{\rm (1)} $X$ is sg-compact.

{\rm (2)} $X$ is a $C_3$-space, i.e., every hsg-closed set is
finite.
\end{lemma}

\begin{lemma}\label{l3}
{\em \cite[Proposition 2.1]{DG1}} For a subset $A$ of a
topological
space $(X,\tau)$ the following conditions are equivalent:

{\rm (1)} $A$ is hsg-closed.

{\rm (2)} $N(X) \cap {\rm int}({\rm cl} (A)) = \emptyset$, where
$N(X)$ denotes the set of nowhere dense singletons in $X$.
\end{lemma}

\section{Sg-compactness in product spaces}\label{s2}

We will start with an example showing that Theorem 2.1 of
\cite{T2} is not true. There, the authors stated (without proof)
that every sg-compact space is go-compact (it is our guess that
they assumed that g-open sets are sg-open).

\begin{example}\label{e1}
{\em Let ${\mathbb N}$ be set of all positive integers. We
consider the following topology $\tau$ on ${\mathbb N}$ given by
$\tau  = \{ \emptyset , {\mathbb N} \} \cup \{ U_n = \{ n, n+1,
n+2, \ldots \} \colon n \geq 3 \}$.

We first show that $({\mathbb N},\tau)$ is sg-compact. Observe
that every singleton of $({\mathbb N},\tau)$ is nowhere dense.
Since every nonempty semi-open set has finite complement,
$({\mathbb N},\tau)$ is semi-compact. By \cite[Remark 2.7
(i)]{DG1}, $({\mathbb N},\tau)$ is sg-compact.

However, every singleton of $({\mathbb N},\tau)$ is g-open, and
so $({\mathbb N},\tau)$ fails to be go-compact.}
\end{example}

\begin{lemma}\label{lmax1}
Let  $X = \prod_{i \in I} X_i$ be a product space. If infinitely
many $X_i$ are not indiscrete, then $X$ contains an infinite
nowhere dense subset.
\end{lemma}

{\em Proof.} Let $J$ be an infinite subset of $I$ such that $X_i$
is not indiscrete for each $i \in J$. We may choose $J$ in such
a way that $I \setminus J$ is also infinite. Then, for each $i
\in J$, there exists a closed set $A_i \subseteq X_i$ distinct
from the empty set and from $X_i$. Now form the product of all
$A_i$, $i \in J$, and of all $X_i$, $i \not\in J$,
and call it $A$. Then $A$ is closed in $X$, infinite and clearly
nowhere dense.  $\Box$

As a consequence of Lemma~\ref{l2} we therefore have:

\begin{corollary}\label{cm2}
If a product space $X = \prod_{i \in I} X_i$ is sg-compact, then
only finitely many $X_i$ are not indiscrete. $\Box$
\end{corollary}

\begin{theorem}\label{t1}
Let $(X_i,\tau_i)_{i \in I}$ be a family of topological spaces.
If the product space $X = \prod_{i \in I} X_i$ is sg-compact,
then either all factor spaces are finite or exactly one of them
is infinite and sg-compact and the rest are finite and locally
indiscrete.
\end{theorem}

{\em Proof.} Suppose that two factor spaces, say $X_i$ and $X_j$,
are infinite. Let $p_i$ denotes the projection from $X$ onto
$X_i$ for any $i \in I$. Let $k \in I$. If $x_k \in X_k$, then
$p_k^{-1} (\{ x_k \} )$ is infinite, hence cannot be nowhere
dense since $X$ is sg-compact. Thus $\{ x_k \}$ is not nowhere
dense in $X_k$. Consequently, each factor space $X_k$ must be
locally indiscrete. By Corollary~\ref{cm2} and Lemma~\ref{l1},
each singleton in $X$ is locally dense and so every subset of $X$
is sg-open. Since $X$ is sg-compact, $X$ must be finite, a
contradiction. Hence, at most one factor space can be infinite.

Now suppose that $X_j$ is infinite and that $X_i$ is finite for
$i \neq j$. For each $x_i \in X_i$, where $i \neq j$, $p_i^{-1}
(\{ x_i \} )$ is infinite, therefore $\{ x_i \}$ cannot be
nowhere dense in $X_i$. So $X_i$ is locally indiscrete for $i
\neq j$. By Corollary~\ref{cm2} and Lemma~\ref{l1} it follows
that for each $x \in X$, $\{ x \}$ is nowhere dense in $X$ if and
only if $\{ x_j \}$ is nowhere dense in $X_j$.

Assume now that $X_j$ is not sg-compact. Then $X_j$ contains an
infinite hsg-closed subset, say $A_j$. Let $A = p_j^{-1} (A_j)$.
We want to show that $N(X) \cap {\rm int}({\rm cl} (A)) =
\emptyset$, where $N(X)$ denotes the set of nowhere dense
singletons in $X$. If there exists a point $x \in N(X) \cap {\rm
int}({\rm cl} (A))$, then $x$ has an open neighbourhood $W$
contained in ${\rm cl} (A)$. Also, $\{ x_j \}$ is nowhere dense
in $X_j$ and $x_j \in p_j(W) \subseteq p_j({\rm cl} (A))
\subseteq {\rm cl} (A_j)$. So $x_j \in {\rm int}({\rm cl}
(A_j))$, a contradiction to the hsg-closedness of $A_j$. Hence,
by Lemma~\ref{l3}, $A$ is hsg-closed and infinite, a
contradiction. Therefore, $X_j$ is sg-compact. $\Box$
\bigskip

Tapi, Thakur and Sonwalkar \cite[Theorem 2.7]{T2} stated our
result for two topological spaces but their proof is wrong as
they claimed the projection mapping being sg-irresolute. They
used the wrong lemma from \cite{T1} that the product of sg-closed
sets is sg-closed. The following example will correct their
claims.

\begin{example}\label{e2}
{\em Let $X = \{ a,b,c \}$ and let $\tau = \{ \emptyset, \{ a,b
\}, X \}$. Set $A = \{ b,c \}$.

(i) First observe that $A$ is sg-closed in $(X,\tau)$ but $A
\times A$ is not sg-closed in $X \times X$, since $A \times A
\subseteq X \times X \setminus \{ (a,c) \}$ and ${\rm scl} (A
\times A) = X \times X$.

(ii) If $p$ is the projection mapping from $X \times X$ onto $X$,
then $p^{-1} (A)$ is not sg-closed in $X \times X$, i.e., the
projection map need not be always sg-irresolute.

(iii) We already noted that if \fxy\ is open and continuous, then
the preimage of every nowhere dense subset of $Y$ is nowhere
dense in $X$. There is no similar result for hsg-closed sets. If
$\sigma$ denotes the indiscrete topology on $X$, then $S =
\{ a,b \}$ is hsg-closed in $(X,\sigma)$ but $q^{-1}(S)$ is not
hsg-closed in $(X,\sigma) \times (X,\tau)$, where $q$ denotes the
projection mapping from $(X,\sigma) \times (X,\tau)$ onto
$(X,\sigma)$.}
\end{example}

The following result shows when the inverse image of a hsg-closed
set is also hsg-closed. Recall that a function \fxy\ is called
{\em almost open} if the image of every regular open set is open.
We say that \fxy\ is {\em anti-$\delta$-open} if the image of
every nowhere dense singleton is nowhere dense. Observe that
if $Y$ is dense-in-itself and $T_D$ (= singletons are locally
dense), then \fxy\ is always anti-$\delta$-open; in particular
every real-valued function is anti-$\delta$-open.

\begin{proposition} 
If \fxy\ is an almost open, continuous, anti-$\delta$-open
surjection, then the inverse image of every hsg-closed set is
hsg-closed.
\end{proposition}

{\em Proof.} Let $B$ be hsg-closed in $Y$ and set $A =
f^{-1}(B)$. If for some nowhere dense singleton $x$ of $X$ we
have $x \in {\rm int} ({\rm cl} (A))$, then $f(x) \in f({\rm int}
({\rm cl} (A))) \subseteq {\rm int} (f({\rm cl} (A))) \subseteq
{\rm int} ({\rm cl} (f(A))) = {\rm int} ({\rm cl} (B))$. Since
$f(x)$ is nowhere dense in $Y$, $B$ is not hsg-closed. By
contradiction, $A$ is hsg-closed. $\Box$

\begin{remark}\label{r1}
{\em (i) Let $A$ be an infinite set with $p \not\in A$. Let $X
= A \cup \{ p \}$ and $\tau = \{ \emptyset, A, X \}$. We
observed in \cite{DG1} that $X \times X$ contains an infinite
nowhere dense subset, so even the finite product of
sg-compact spaces need not be sg-compact.

(ii) It is rather unexpected that the projection map fails to
be sg-irresolute in general, since it is always irresolute and
gs-irresolute.}
\end{remark}

The two examples of infinite sg-compact spaces in \cite{DG1} and
the infinite sg-compact space from Example~\ref{e1} are not even
weakly Hausdorff (however one of them is $T_1$). As every
hereditarily compact kc-space must be finite, it is natural to
ask whether there are any infinite sg-compact semi-Hausdorff
spaces (there do exist infinite hereditarily compact
semi-Hausdorff spaces). Recall here that a topological space
$(X,\tau)$ is called {\em semi-Hausdorff} \cite{MP1} if every two
distinct points of $X$ can be separated by disjoint semi-open
sets.

Recall additionally that a space $(X,\tau)$ is called {\em
hyperconnected} if every open subset of $X$ is dense, or
equivalently, every pair of nonempty open sets has nonempty
intersection. In the opposite case $X$ is called {\em
hyperdisconnected}. If every infinite open subspace of $X$ is
hyperdisconnected, then we will say that $X$ is {\em
quasi-hyperdisconnected}. Note that not only Hausdorff spaces but
also semi-Hausdorff spaces are quasi-hyperdisconnected (but not
vice versa).

\begin{proposition}\label{p2}
Every quasi-hyperdisconnected sg-compact space $(X,\tau)$ is
finite.
\end{proposition}

{\em Proof.} Assume that $X$ is infinite. Let $U$ and $V$ be
disjoint non-empty open subsets of $X$. Note that
either $X \setminus U$ or $X \setminus V$ is infinite. Assume
that $X \setminus U$ is infinite. Since ${\rm cl} (U) \setminus
U$ is hsg-closed (in fact even nowhere dense), by Lemma~\ref{l2},
${\rm cl} (U) \setminus U$ is finite and hence $X \setminus {\rm
cl} (U)$ is infinite and open. Set $A_1 = U$. Since $X$ is
quasi-hyperdisconnected, proceeding as above, we can construct
an open subset of $X \setminus {\rm cl} (U)$ and hence of $X$,
say $U_2$, such that the complement of the closure of $U_2$ in
$X \setminus {\rm cl} (A_1)$ is infinite. Using the method above,
we can construct an infinite pairwise disjoint family $A_1, A_2,
\ldots$ of non-empty open subsets of $(X,\tau)$. Since sg-compact
spaces are semi-compact and thus satisfy the finite chain
condition,
$X$ must be finite.  $\Box$

\begin{corollary}
Every sg-compact, semi-Hausdorff space is finite.
\end{corollary}

We have just seen that under some very low separation axioms,
sg-compact spaces very easily become finite. If we replace the
weak separation axiom with a weaker form of strong
irresolvability, we again have finiteness. By definition, a
nonempty topological space $(X,\tau)$ is called {\em resolvable}
\cite{H1} if $X$ is the disjoint union of two dense (or
equivalently codense) subsets. In the opposite case $X$ is called
{\em irresolvable}. A topological space $(X,\tau)$ is {\em
strongly irresolvable} \cite{FL1} if no nonempty open set is
resolvable.

\begin{proposition}\label{p3}
Every sg-compact space $(X,\tau)$ which is the topological sum
of a locally indiscrete space and a strongly irresolvable space
is finite.
\end{proposition}

{\em Proof.} We will use a result in \cite{G1} which states that
a space is finite if and only if every cover by $\beta$-open sets
(i.e., sets which are dense in some regular closed subspace) has
a finite subcover. If $\cal U$ is a cover of $X$ by $\beta$-open
sets, then by \cite[Theorem 2.1]{CGR1} every element of $\cal U$
is sg-open. Since $X$ is sg-compact, $\cal U$ has a finite
subcover. This shows that $X$ is finite. $\Box$
\bigskip

We already mentioned in Remark~\ref{r1} that the product of two
sg-compact spaces need not be sg-compact. Thus we have the
natural question: When is the product of two sg-compact spaces
also sg-compact? What turns out is that only in one very special
case the product of a sg-compact space with another sg-compact
space is also sg-compact. First we note a result whose proof is
easy and hence omitted.

\begin{proposition}\label{pl1}
Let $(X_{\alpha},{\tau}_{\alpha})_{\alpha \in \Omega}$ be a
family of topological spaces. For the topological sum $X =
\sum_{\alpha \in \Omega} X_{\alpha}$ the following conditions are
equivalent:

{\rm (1)} $X$ is a sg-compact space.

{\rm (2)} Each $X_{\alpha}$ is a sg-compact space and $|\Omega|
< \aleph_{0}$.
\end{proposition}

\begin{lemma}\label{lem}
Let $(X,\tau)$ be any space and let $(Y,\sigma)$ be indiscrete.
Let $A \subseteq X \times Y$ and let $p: X \times Y \rightarrow
X$ denote the projection. Then ${\rm int} ({\rm cl} (A)) = {\rm
int} ({\rm cl} (p(A))) \times Y$.
\end{lemma}

{\em Proof.} If $(x,y) \in {\rm int} ({\rm cl} (A))$, there
exists an open neighbourhood $U_x$ of $x$ such that $U_x \times
Y \subseteq {\rm cl} (A)$. Then $x \in p(U_x) \subseteq {\rm cl}
(p(A))$ and so $(x,y) \in {\rm int} ({\rm cl} (p(A))) \times Y$.

Now, let $x \in {\rm int} ({\rm cl} (p(A)))$ and $y \in Y$.
Choose an open set $U_x \subseteq X$ containing $x$ such that
$U_x \subseteq {\rm cl} (p(A))$. We claim that $U_x \times Y
\subseteq {\rm cl} (A)$. Suppose there is a point $(x',y') \in
U_x \times Y$ not in ${\rm cl} (A)$. Then there exists an open 
set $W_{x'} \subseteq U_x$ containing $x'$ such that $(W_{x'}
\times Y) \cap A = \emptyset$. Consequently, $W_{x'} \cap p(A)
= \emptyset$, a contradiction. Hence, $(x,y) \in {\rm int} ({\rm
cl} (A))$. $\Box$

\begin{theorem}\label{t2}
If $(X,\tau)$ is sg-compact and $(Y,\sigma)$ is finite and
locally indiscrete, then $X \times Y$ is sg-compact.
\end{theorem}

{\em Proof.} Since $Y$ is a finite topological sum of indiscrete
spaces, by Proposition~\ref{pl1} it suffices to assume that $Y$
is indiscrete. Suppose that $A \subseteq X \times Y$ is infinite
and hsg-closed. Then, $p(A)$ is infinite and hence, by
Lemma~\ref{l2} and Lemma~\ref{l3}, we have $N(X) \cap {\rm int}
({\rm cl} (p(A))) \neq \emptyset$. Pick $x \in N(X) \cap {\rm
int} ({\rm cl} (p(A)))$ and $y \in Y$. Then, $\{ (x,y) \}$ is
nowhere dense in $X \times Y$ and, by Lemma~\ref{lem}, we have
$(x,y) \in {\rm int} ({\rm cl} (A))$, a contradiction to the
hsg-closedness of $A$. Thus $X \times Y$ is sg-compact. $\Box$

\
E-mail: {\tt dontchev@cc.helsinki.fi}, {\tt
ganster@weyl.math.tu-graz.ac.at}
\
\end{document}